\documentclass[reqno]{amsart}
\usepackage{amssymb}
\usepackage{amsmath,amscd}
\usepackage{graphicx}
\usepackage{color,hyperref}
\usepackage{tikz}
\usetikzlibrary{matrix,arrows}
\usepackage{enumitem}

\newcommand{\mb}{\mathbf}

\newcommand{\C}{\mathbb C}

\numberwithin{equation}{section}

\theoremstyle{plain}
\newtheorem{theorem}{Theorem}

\theoremstyle{definition}

\newtheorem{example}[theorem]{Example}

\begin{document}
	
\title[Generalized Chebyshev Acceleration]{Generalized Chebyshev Acceleration}

\author[Nurgül Gökgöz]{Nurgül Gökgöz}
\address{Çankaya University, Mathematics Department, 06815 Ankara, Turkey.}
\email{ngokgoz@cankaya.edu.tr}

\date{\today}
	
\begin{abstract}
We use generalized Chebyshev polynomials, associated with the root system $A_2$, to provide a new semi-iterative method for accelerating simple iterative methods for solving linear systems. We apply this semi-iterative method to the Jacobi method, and give an example. There are certain restrictions but the resulting acceleration is rather high.
\end{abstract}
	
	
	
	
\maketitle
Let $A$ be an $n\times n$ invertible matrix with real entries. Consider the linear system $A\mb{x}=\mb{b}$. For a general matrix with reasonable size and entries, we normally use some direct method of solution, such as Gaussian elimination. However, for larger matrices, we may prefer to use an iterative method. Iterative methods can be useful for linear problems involving many variables, where direct methods are computationally expensive, if not impossible. See \cite{greenbaum}, and \cite{varga}, for several iterative algorithms.
	
The Chebyshev semi-iterative method is an enhancement of simple iterative methods for solving linear systems. It improves convergence by utilizing Chebyshev polynomials to optimize the choice of iteration parameters, thereby accelerating convergence compared to standard stationary iterative methods. The method requires an estimate of the smallest and largest eigenvalues of the iteration matrix to construct the Chebyshev polynomial sequence \cite{foxparker}, \cite{varga}.
	
In this paper, we use generalized Chebyshev polynomials associated with the root system $A_2$, see \cite{munthe}, and provide a new approach for accelerating convergence. We apply this new approach to the Jacobi method. Even though, there are certain restrictions, the resulting acceleration is rather high
	
The organization of the paper is as follows. In the first section, we summarize the basics of semi-iterative methods, and fix some notation that will be used in the rest of the paper. In the second section, we outline the classical Chebyshev acceleration method to make its generalization more accessible. In the third section, we provide certain properties of generalized Chebyshev polynomials associated with the root system $A_2$. In the last section, we present our main result for accelerating the convergence of simple iterative methods by using the generalized Chebyshev polynomials. 
	
\section{Semi-Iterative Methods}
An iterative method is a mathematical procedure that uses an initial value to generate a sequence of improving approximate solutions. This general approach can be used to study the solution of the linear systems $A\mb{x}=\mb{b}$. For a comprehensive treatment of the subject, one can see \cite{greenbaum}, and \cite{varga}. In this section, we summarize the basics of semi-iterative methods, and fix some notation that will be used in the rest of the paper. We will mostly follow the exposition of \cite{varga}.
	
Before we can explain semi-iterative methods, we first need an iterative method. For simplicity, we consider the Jacobi method. Given a system $A\mb{x}=\mb{b}$, with nonzero diagonal, we set $A=L+D+U$ where $L$ and $U$ are respectively strictly lower and upper triangular matrices, and $D$ is the diagonal matrix. Consider $M=-D^{-1}(L+U)$ and $\mb{g}=D^{-1}\mb{b}$. The Jacobi method uses the following iterative process
\begin{equation}\label{iteration}
\mb{x}^{(m)}=M\mb{x}^{(m-1)}+\mb{g}, \quad m \geq 1.
\end{equation}
To investigate the convergence, we observe that the required solution of the system $A\mb{x}=\mb{b}$ shall satisfy
\begin{equation}\label{relation}
\mb{x}=M\mb{x}+\mb{g}.
\end{equation}
We define the error vector as $\epsilon^{(m)}:=\mb{x}-\mb{x}^{(m)}$. Combining (\ref{iteration}) and (\ref{relation}), we see that the error vector has the following property
\begin{equation}\label{error}
	\epsilon^{(m)}= M\epsilon^{(m-1)}  =\ldots  = M^m\epsilon^{(0)}.
\end{equation}
The error vector converges to zero if and only if all the eigenvalues $M$ have absolute values less than one, i.e. the spectral radius $\rho(M)$ is less than one. This is easily proved in the generic case. Suppose that $\lambda_1, \ldots \lambda_n$ are eigenvalues and $\{\mb{v}_1,\ldots \mb{v}_n\}$ is a full set of independent eigenvectors of $M$, respectively. We have
\begin{equation}\label{eigenvector}
	\epsilon^{(0)} = \sum_{j=1}^n\alpha_j\mb{v}_j, \quad
	\epsilon^{(m)}=M^m\epsilon^{(0)} = \sum_{j=1}^n\alpha_j\lambda_j^m\mb{v}_j.
\end{equation}
The coefficients $\alpha_j$ depend on the initial guess $\mb{x}^{(0)}$, and therefore they are random. However, the convergence is not related with these coefficients. More precisely, we have $M^m\epsilon^{(0)} \rightarrow 0$ if and only if $|\lambda_j|<1$ for all $j$. 
	
If the absolute value of an eigenvalue, say $\lambda_1$, is greater than the absolute values of all other eigenvalues, then the dominant term in $M^m\epsilon^{(0)}$ is $\alpha_1\lambda_1^m\mb{v}_1$. The rate of convergence ultimately depends on the size of $|\lambda_1|$, and therefore the spectral radius.
	
We may obtain a faster rate of convergence by weighting the components of the error vector, see (\ref{eigenvector}), suitably. Consider the following linear combination
\begin{equation}\label{y}
	\mb{y}^{(m)} := \sum_{j=0}^{m} \nu_j(m)\mb{x}^{(j)}, \quad m\geq 0.
\end{equation}
This is called a \textit{semi-iterative method} with respect to the iterative method of (\ref{iteration}). We want this new iteration to converge to the same solution, so there is an obvious restriction on the coefficients $\nu_j(m)$. If the first guess for the original iteration is correct, i.e. $\mb{x}^{(0)} = \mb{x}$, then $\mb{x}^{(1)}=\mb{x}^{(2)}=\ldots=\mb{x}$, as well. Since we want $\mb{y}^{(m)}$ to converge to $\mb{x}$, we must have
\begin{equation}\label{restrict}
	\sum_{j=0}^{m} \nu_j(m)=1.
\end{equation}
We shall relate the error of the former iterative process with the error of this new setup. For this purpose, we introduce a new notation, and define the error vector $\eta^{(m)}:=\mb{x}-\mb{y}^{(m)}$ in the semi-iterative case. Using (\ref{y}), and (\ref{restrict}), we obtain 
\begin{equation*}
	\eta^{(m)} =\mb{x}-\sum_{j=0}^{m} \nu_j(m)\mb{x}^{(j)} = \sum_{j=0}^{m} \nu_j(m)(\mb{x}-\mb{x}^{(j)}).
\end{equation*}
Recalling $\epsilon^{(j)}=\mb{x}-\mb{x}^{(j)}$, and applying (\ref{error}), we get 
\begin{equation*}
	\eta^{(m)}=\sum_{j=0}^{m} \nu_j(m)M^j\epsilon^{(0)}.	
\end{equation*}
It is now apparent that $\eta^{(m)}$ can be obtained by applying a suitable polynomial of $M$ to the initial error vector $\epsilon^{(0)}$ of the iterative method of (\ref{iteration}). We define
\begin{equation*}
	p_m(\lambda) := \sum_{j=0}^{m} \nu_j(m)\lambda^j.
\end{equation*}
Using the setting as in (\ref{eigenvector}), we easily see that
\begin{equation}\label{eta}
	\eta^{(m)} = \sum_{j=1}^n\alpha_j p_m(\lambda_j)\mb{v}_j = p_m(M)\epsilon^{(0)}.
\end{equation}
With standard definitions of vector and matrix norms, see \cite{varga}, it follows that 
\begin{equation*}
	\lVert\eta^{(m)}\rVert = \lVert  p_m(M)\epsilon^{(0)} \rVert \leq \lVert  p_m(M)\rVert \cdot \lVert\epsilon^{(0)} \rVert, \quad m\geq 0.
\end{equation*}
If we can make the quantities $\lVert p_m(M)\rVert$ smaller, compared to $\lVert M^m\rVert$, then the bounds for $\lVert\eta^{(m)} \rVert$ for the semi-iterative method are proportionally smaller, compared to $\lVert\epsilon^{(m)}\rVert$.
	
\section{Chebyshev Acceleration}
In this section, we give a short summary of the Chebyshev semi-iterative method with respect to (\ref{iteration}) following \cite{varga}. In the following sections, we will generalize this construction to the generalized Chebyshev polynomials associated with the root system $A_2$.
	
Suppose that $M$ has real eigenvalues with absolute values strictly less than one. Suppose also that we know the spectral radius $0<\rho=\rho(M)<1$. By the definition of Chebyshev semi-iterative method, we choose
\begin{equation}\label{pm}
	p_m(\lambda)=\frac{C_m(\lambda/\rho)}{C_m(1/\rho)}.
\end{equation}
where $C_m$ is the $m$th Chebyshev polynomial. We obviously have $p_m(1)=1$. In other words the restriction (\ref{restrict}) is satisfied. The Chebyshev polynomials $C_m$ are uniquely defined by the following functional equation
	\[C_m(\cos(\theta))=\cos(m\theta).\]
The first few Chebyshev polynomials are $C_0(t)=1$, $C_1(t)=t$, $C_2(t)=2t^2-1$, and they satisfy the recurrence relation 
	\[C_{m}(t)=2tC_{m-1}(t)-C_{m-2}(t), \quad m\geq 2.\]
This relation, applied to $C_m(\lambda/\rho) = C_m(1/\rho)p_m(\lambda)$, gives us
	\[p_m(\lambda)C_m(1/\rho) = 2\frac{\lambda}{\rho} C_{m-1}(1/\rho)p_{m-1}(\lambda) -C_{m-2}(1/\rho)p_{m-2}(\lambda), \quad m\geq 2.\] 
We observe that
	\[p_m(\lambda) = \frac{2\lambda C_{m-1}(1/\rho)}{C_m(1/\rho)}p_{m-1}(\lambda) - \frac{C_{m-2}(1/\rho)}{C_m(1/\rho)}p_{m-2}(\lambda), \quad m \geq 2.\]
We now plug in $M$ for $\lambda$, and postmultiply throughout by $\epsilon^{(0)}$. Using the relation with the error vector, see (\ref{eta}), we obtain 
	\[\eta^{(m)} = \frac{2C_{m-1}(1/\rho)}{C_m(1/\rho)}M\eta^{(m-1)} - \frac{C_{m-2}(1/\rho)}{C_m(1/\rho)}\eta^{(m-2)}, \quad m \geq 2.\]
Recall that $\eta^{(j)}=\mb{x}-\mb{y}^{(j)}$ for each $j$, and $M\mb{x}+\mb{g}=\mb{x}$. Thus, we obtain a new iterative process,
\begin{equation}\label{itercheby}
	\mb{y}^{(m)} = \frac{2C_{m-1}(1/\rho)}{\rho C_{m}(1/\rho)}(M\mb{y}^{(m-1)}+\mb{g}) - \frac{C_{m-2}(1/\rho)}{C_{m}(1/\rho)}\mb{y}^{(m-2)}.
\end{equation}
For simplicity, we can choose $\mb{y}^{(0)}=\mb{x}^{(0)}$, and  $\mb{y}^{(1)}=\mb{x}^{(1)}$. This is of a similar form to the original iteration, see (\ref{iteration}), but involves two earlier vectors instead of one. We also note that, the terms $\mb{y}^{(j)}$ for $j=0,1,2, \ldots$ can be computed directly without going through (\ref{iteration}). On the other hand, this semi-iterative method requires an additional vector storage. This is of minor importance in computer applications. 
	
\begin{example}\label{ex1}
As an example, consider the following system
\[\left[\begin{array}{cccc}
	576&0&0&1\\
	144&144&0&5\\
	0&144&144&25\\
	0&0&1&1
\end{array}\right] \left[\begin{array}{c}
			x_1\\x_2\\x_3\\x_4
\end{array}\right]=\left[\begin{array}{c}
			577\\293\\313\\2
\end{array}\right].\]
Note that this system has the unique solution $\mb{x}=(1,1,1,1)$. According to the iterative method of Jacobi, we have
\[M=\left[\begin{array}{cccc}
	0&0&0&-1/576\\
	-1&0&0&-5/144\\
	0&-1&0&-25/144\\
	0&0&-1&0
	\end{array}\right]\quad\text{and}\quad \mb{g}=\left[\begin{array}{c}
			577/576 \\ 293/144\\ 313/144\\ 2
\end{array}\right]. \]
The eigenvalues of the iteration matrix $M$ are $-1/2$, $1/4$, $1/6$ and $1/12$. We have $\rho=1/2<1$. The following table gives the first few iterates and the norm of the related errors under (\ref{iteration}) starting from the zero vector. The error in $\mb{x}^{(m)}$ is ultimately decreased at each stage by the factor $\rho=1/2$.  
\[\begin{array}{|c|cccc|c|} \hline
m&&\mb{x}^{(m)}& & &\lVert\epsilon^{(m)}\rVert \\ \hline
	0 &  0 & 0 & 0 & 0   & 2.000  \\
	1 &  1.001 & 2.034 & 2.173 & 2.000   & 1.856  \\
	2 &  0.998 & 0.963 & -0.208 & -0.173   & 1.684  \\
	3 &  1.002 & 1.042 & 1.240 & 2.208   & 1.232  \\
	4 &  0.997 & 0.956 & 0.747 & 0.759   & 0.351 \\
	5 &  1.000 & 1.010 & 1.085 & 1.252   & 0.266 \\
	6 &  0.999 & 0.990 & 0.945 & 0.914   & 0.101 \\
	7 &  1.000 & 1.003 & 1.024 & 1.054   & 0.059 \\
	8 &  0.999 & 0.997 & 0.987 & 0.975   & 0.027 \\ \hline
\end{array}\]
		
The ultimate behavior of the error in $\mb{y}^{(m)}$ requires a little more consideration. The following table gives the first few iterates, and the norm of the related errors, under the semi-iterative process (\ref{itercheby}) starting from the zero vector.
\[\begin{array}{|c|cccc|c|}\hline
	m& &\mb{y}^{(m)}&&&\lVert\eta^{(m)}\rVert\\ \hline
	0 &    0 & 0 & 0 & 0   & 2.000 \\
	1 &    1.001 & 2.034 & 2.173 & 2.000   & 1.856 \\
	2 &    1.140 & 1.101 & -0.238 & -0.198   & 1.731 \\
	3 &    1.002 & 0.813 & 1.024 & 2.256   & 1.270 \\
	4 &    0.987 & 0.943 & 1.055 & 1.059   & 0.099 \\
	5 &    0.999 & 1.024 & 1.047 & 0.850   & 0.158\\
	6 &    1.001 & 1.009 & 0.997 & 0.944   & 0.056\\
	7 &    1.000 & 0.999 & 0.996 & 1.013   & 0.013\\
	8 &    0.999 & 0.998 & 0.998 & 1.007   & 0.008 \\ \hline
\end{array}\]
To analyze the ultimate behavior of the error, we follow the treatment given in \cite[7.6]{foxparker}. Let $\theta$ be the positive real number such that $1/\rho=\cosh \theta$. In our case, we have $\theta \approx 1.317$. Expressing (\ref{itercheby}) in terms of the hyperbolic cosine function, and focusing on the dominant terms, we see that this iterative process is ultimately equivalent to the following 
		\begin{equation*}
			\mb{y}^{(m)} = (1+e^{-2\theta})(M\mb{y}^{(m-1)}+\mb{g}) - e^{-2\theta}\mb{y}^{(m-2)}.
		\end{equation*}
		Moreover, the error vector satisfies the iteration 
		\begin{equation*}
			\eta^{(m)} = (1+e^{-2\theta})M\eta^{(m-1)} - e^{-2\theta}\eta^{(m-2)},
		\end{equation*}
		which can be put into the matrix form as follows
		\[\left[\begin{array}{c}
			\eta^{(m)}\\\eta^{(m-1)}
		\end{array}\right] = \left[\begin{array}{cc}
			(1+e^{-2\theta})M&e^{-2\theta}\\\mb{I}&\mb{0}
		\end{array}\right] \left[\begin{array}{c}
			\eta^{(m-1)}\\\eta^{(m-2)}
		\end{array}\right].\]
The ultimate convergence depends on the matrix on the right-hand side. It is easy to see that the eigenvalues $\mu$ of this matrix satisfy the equation
	\[\mu^2-\mu(1+e^{-2\theta})\lambda + e^{-2\theta}=0,\]
where $\lambda$ is an eigenvalue of $M$. It turns out that the largest $|\mu|$ corresponds to the largest $|\lambda|=\rho$. After some manipulation, one can deduce
	\[|\mu|_\text{max} = \frac{1-\sqrt{1-\rho^2}}{\rho}.\]
In our example, this quantity equals $e^{-\theta} \approx 0.268$, so the iterative process has improved very significantly. Due to the non-Hermitian nature of the example, the ratios $\lVert\eta^{(m+1)}\rVert / \lVert\eta^{(m)}\rVert$ do not behave steadily. However, we have experimentally verified their geometric mean approaches to this number.
\end{example}
	
\section{Generalized Chebyshev Polynomials}
The generalized Chebyshev polynomials are defined by means of exponential invariants of Bourbaki \cite{Bourbaki}. These polynomials are studied by Veselov \cite{veselov}, and by Hoffman and Withers \cite{hoffwith}, independently, and they have astonishing properties in the theory of dynamical systems and numerical analysis. 
	
\begin{figure}[htbp]
		\centering
		\includegraphics[scale=0.7]{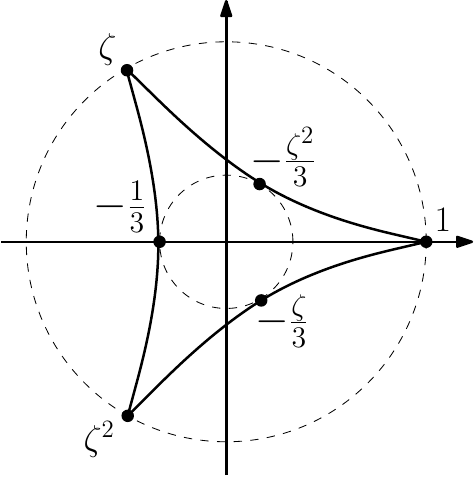}
		\caption{The region $\Delta$.}
		\label{fig:deltoid}
	\end{figure}

In this section, we provide certain properties of the generalized Chebyshev polynomials associated with the root system $A_2$. We will use a normalized version given in \cite{munthe}. For $\theta=(\theta_1,\theta_2) \in \mathbb{C}^2$, we set 
	\begin{equation}\label{gencos}
		\varphi_1(\theta)=\frac{1}{3}\left(e^{2\pi i \theta_1}+e^{-2\pi i \theta_2}+e^{2\pi i (\theta_2-\theta_1)}\right), \quad \varphi_2(\theta)=\overline{\varphi_1(\theta)}.
	\end{equation}
The function $\Phi(\theta)=(\varphi_1(\theta),\varphi_2(\theta))$ is called \textit{the generalized cosine} associated with the root system $A_2$. Similar to the univariate Chebyshev polynomials, and their relation to the cosine function, there exists a family of bivariate polynomial mappings $\mathcal{A}_m$ that satisfies the following functional identity:
	\begin{equation*}
		\mathcal{A}_m(\Phi(\theta)) = \Phi(m\theta).
	\end{equation*}
If we write $\mathcal{A}_m=(f_m,g_m)$, then the first component satisfies $f_m(\Phi(\theta))=\varphi_1(m\theta)$. Putting $x=\varphi_1(\theta)$, and $\bar{x}=\varphi_2(\theta)$, and using (\ref{gencos}), we see that the second component is obtained by switching the variables of the first component, i.e. $f_m(x,\bar{x})=g_m(\bar{x},x)$. In this paper, we will use the first component, and we write $f_m(x)$ instead of $f_m(x,\bar{x})$, for simplicity. The first few of these polynomials are listed below:
	\begin{align*}
		f_0(x) &= 1\\
		f_1(x) &= x\\
		f_2(x) &= 3x^2-2\bar{x}\\
		f_3(x) &= 9x^2+9x\bar{x}+1
	\end{align*}
	These bivariate polynomials satisfy the following recursive relation
	\begin{equation}\label{recursion}
		f_m = 3xf_{m-1}-3\bar{x}f_{m-2}+f_{m-3}, \quad m\geq 3.
	\end{equation}
The key property of these polynomials that motivates this work is the fact that $f_m$ maps a large subset of $\C$ onto itself. Consider
	\begin{equation}\label{deltoid}
		\Delta=\{ \varphi_1(\theta) \mathrel{|} \theta \in \mathbb{R}^2\}.
	\end{equation}
This region is enclosed by \textit{Steiner hypocycloid}, a deltoid curve with three corners at third roots of unity. For an illustration of this region, see Figure~\ref{fig:deltoid}. One can alternatively write that
	\[\Delta=\{ x+iy \mathrel{|} 3\left(x^{2}+y^{2}+1\right)^{2}+8\left(-x^{3}+3xy^{2}\right) \leq 4\}.\]
	
It is a useful fact that  $\Delta\times \Delta$ is the Julia set of the polynomials $\mathcal{A}_m=(f_m,g_m)$. See \cite{veselov} or \cite{veselov-survey}. So the region $\Delta$ is analogous to the interval $[-1,1]$ in the single variable case. If $\gamma$ is a real number in the interval $[-1,1]$, then the sequence
	\[C_1(\gamma),C_2(\gamma),C_3(\gamma),\ldots\]
remains in the interval $[-1,1]$. Similarly, if we have a complex number $\tau\in\Delta$, then the sequence 
	\[f_1(\tau),f_2(\tau),f_3(\tau),\ldots\]
remains in the region $\Delta$. This analogy allows us to generalize the Chebyshev semi-iterative method to a wider setting. Following (\ref{pm}), and adapting to this new situation, we choose
	\[p_m(\lambda) = \frac{f_m(\lambda/\lambda_1)}{f_m(1/\lambda_1)},\]
where $\lambda_1$ is an eigenvalue of $M$ with maximal absolute value. We obviously have $p_m(1)=1$, and the basic restriction (\ref{restrict}) is satisfied. Moreover, using the recursive relation (\ref{recursion}), we now have
	\begin{align}\label{keyrecursion}
		\begin{split}	
			f_m(1/\lambda_1)p_m(\lambda) =\ & 3 \frac{\lambda}{\lambda_1} f_{m-1}(1/\lambda_1)p_{m-1}(\lambda)\\
			&-3\frac{\bar{\lambda}}{\bar{\lambda}_1} f_{m-2}(1/\lambda_1)p_{m-2}(\lambda)\\
			&+f_{m-3}(1/\lambda_1)p_{m-3}(\lambda), \quad m\geq 3.
		\end{split}	
	\end{align}
We want the condition (\ref{eta}) to be satisfied. In order to relate the term $\bar{\lambda}$ in the above equation to the matrix $M$, while forming our new interative method, we need to introduce a new matrix that is closely related to $M$. This is described in the next section. 
	
\section{Main Result}
In this final section, we present our main result, \textit{a generalized Chebyshev semi-iterative method} with respect to (\ref{iteration}). 
	
Suppose that $M$ has spectral radius $0<\rho(M)<1$ with possibly non-real eigenvalues. Suppose that we know an eigenvalue $\lambda_1$ with maximal absolute value, i.e. $0<|\lambda_1|=\rho(M)$. Suppose further that $\lambda/\lambda_1$ lies in the region $\Delta$, defined by (\ref{deltoid}), for each eigenvalue $\lambda$ of $M$.
	
Suppose $\widetilde{M}$ is an $n\times n$ matrix with the same set of eigenvalues as $M$ but the eigenvectors are switched with their complex conjugate counterparts. If $M$ is a normal matrix then $\widetilde{M}$ can be chosen to be the conjugate transpose of $M$. More generally, if $D_0=P^{-1}MP$ is a diagonal matrix, then one can pick $\widetilde{M} = P \overline{D}_0 P^{-1}$. We also suppose that $\widetilde{M}\mb{x} + \widetilde{\mb{g}}=\mb{x}$ for some vector $\widetilde{\mb{g}}$. This condition allows us to deal with complex conjugation occurring in (\ref{keyrecursion}), and to have (\ref{eta}) to be satisfied for this particular $p_m$. 
	
Let us plug in $M$, and $\widetilde{M}$ for $\lambda$, and $\bar{\lambda}$, respectively, in (\ref{keyrecursion}), and postmultiply throughout by $\epsilon^{(0)}$. Using the relation with the error vector, see (\ref{eta}), and by the choice of $\widetilde{M}$ associated with $M$, we obtain 
	\begin{align*}
		\eta^{(m)}=\ &\frac{3f_{m-1}(1/\lambda_1)}{\lambda_1 f_m(1/\lambda_1)}M\eta^{(m-1)} \\
		& -\frac{3f_{m-2}(1/\lambda_1)}{\bar{\lambda}_1 f_m(1/\lambda_1)} \widetilde{M}\eta^{(m-2)}\\
		& +\frac{f_{m-3}(1/\lambda_1)}{f_m(1/\lambda_1)}\eta^{(m-3)}, \quad m\geq 3.
	\end{align*}
Recall that $\eta^{(j)}=\mb{x}-\mb{y}^{(j)}$ for each $j$. With standard definitions of vector and matrix norms, see \cite{varga}, it follows that  
	\begin{equation*}
		\eta^{(m)} = \lVert  p_m(M)\epsilon^{(0)} \rVert \leq \lVert  p_m(M)\rVert \cdot \lVert\epsilon^{(0)} \rVert, \quad m\geq 0.
	\end{equation*}
Note that the term $|1/f_m(1/\lambda_1)|$ decreases rapidly, and so does $\lVert p_m(M)\rVert$. Therefore, this new method has a big potential for accelerating simple iterative methods. 
	
We already have $M\mb{x}+\mb{g}=\mb{x}$ in the original setup. In addition, we have assumed  $\widetilde{M}\mb{x} + \widetilde{\mb{g}}=\mb{x}$ for some vector $\widetilde{\mb{g}}$. Thus, we obtain the following iterative method,
	\begin{align}\label{main}
		\begin{split}		
			\mb{y}^{(m)}=\ &\frac{3f_{m-1}(1/\lambda_1)}{\lambda_1 f_m(1/\lambda_1)}\left( M\mb{y}^{(m-1)} + \mb{g} \right) \\
			& -\frac{3f_{m-2}(1/\lambda_1)}{\bar{\lambda}_1 f_m(1/\lambda_1)}\left( \widetilde{M}\mb{y}^{(m-2)} + \widetilde{\mb{g}} \right)\\
			& +\frac{f_{m-3}(1/\lambda_1)}{f_m(1/\lambda_1)}\mb{y}^{(m-3)}, \quad m\geq 3.
		\end{split}		
	\end{align}
For simplicity, we choose $\mb{y}^{(0)}=\mb{x}^{(0)}$, $\mb{y}^{(1)}=\mb{x}^{(1)}$, and $\mb{y}^{(2)}=\mb{x}^{(2)}$. This new semi-iterative method is of a similar form to (\ref{itercheby}), but involves three earlier vectors instead of two. This requires an additional vector storage.

\begin{example}\label{ex2}
		As an example, consider the following system
		\[\left[\begin{array}{cccc}
			2250&0&0&17\\                                                                                  2250&2250&0&181\\                                                                              0&900&900&53\\                                                                                 0&0&1&1
		\end{array}\right] \left[\begin{array}{c}
			x_1\\x_2\\x_3\\x_4
		\end{array}\right]=\left[\begin{array}{c}
			2267\\4681\\1853\\2
		\end{array}\right].\]
Note that this system has the unique solution $\mb{x}=(1,1,1,1)$. According to the iterative method of Jacobi, we have
		\[M=\left[\begin{array}{cccc}
			0&0&0&-17/2250\\
			-1&0&0&-181/2250\\
			0&-1&0&-53/900\\
			0&0&-1&0
		\end{array}\right]\quad\text{and}\quad \mb{g}=\left[\begin{array}{c}
			2267/2250\\ 4681/2250\\ 1853/900\\ 2
		\end{array}\right]. \]
The eigenvalues of the iteration matrix $M$ are $-1/2$, $1/10$, and $1/5\pm i/3$. We have $\rho=1/2<1$. The following table gives the first few iterates and the norm of the related errors under (\ref{iteration}) starting from the zero vector. The error in $\mb{x}^{(m)}$ is ultimately decreased at each stage by the factor $\rho=1/2$.  
\[\begin{array}{|c|cccc|c|} \hline
	m&&\mb{x}^{(m)}& & &\lVert\epsilon^{(m)}\rVert \\ \hline
	0 &  0 & 0 & 0 & 0   & 2.000 \\
	1 &  1.007 & 2.080 & 2.058 & 2.000   & 1.813\\
	2 &  0.992 & 0.912 & -0.139 & -0.058 & 1.557  \\
	3 &  1.008 & 1.092 & 1.150 & 2.139   & 1.152  \\
	4 &  0.991 & 0.900 & 0.840 & 0.849   & 0.241  \\
	5 &  1.001 & 1.020 & 1.108 & 1.159   & 0.194  \\
	6 &  0.998 & 0.986 & 0.969 & 0.891   & 0.113   \\
	7 &  1.000 & 1.009 & 1.020 & 1.030   & 0.037  \\
	8 &  0.999 & 0.996 & 0.988 & 0.979   & 0.023\\ \hline
\end{array}\]

The computation of $\mb{y}^{(m)}$ requires some additional information. First, note that if we divide an eigenvalue $\lambda$ of $M$ by $\lambda_1 = -1/2$, the resulting value $\lambda/\lambda_1$ lies in the region $\Delta$. Recall that $\widetilde{M}$ is an $n\times n$ matrix with the same set of eigenvalues as $M$ but the eigenvectors are switched with their complex conjugate counterparts.
We pick $\widetilde{M} = P \overline{D}_0 P^{-1}$, where $P$ is obtained by the distinct eigenvectors of $M$, and $D_0=P^{-1}MP$ is a diagonal matrix whose diagonal entries are the eigenvalues of $M$. We also pick $\widetilde{\mb{g}}$, a vector that satisfies the relation $\widetilde{M}\mb{x} + \widetilde{\mb{g}}=\mb{x}$. 
		
The following table gives the first few iterates, and the norm of the related errors, under the generalized Chebyshev semi-iterative method (\ref{main}).
\[\begin{array}{|c|cccc|c|}\hline
	m& &\mb{y}^{(m)}&&&\lVert\eta^{(m)}\rVert\\ \hline
	0&0 & 0 & 0 & 0   & 2.000 \\
	1&1.007 & 2.080 & 2.0588 & 2.000   & 1.813 \\ 
	2&0.992 & 0.912 & -0.139 & -0.058   & 1.557 \\
	3&1.013 & 1.074 & 1.118 & 1.758   & 0.771  \\
	4&0.997 & 0.960 & 0.933 & 0.924   & 0.108 \\
	5&1.000 & 1.004 & 1.019 & 1.031   & 0.037 \\
	6&0.999 & 0.998 & 0.997 & 0.992   & 0.008 \\
	7&1.000 & 1.000 & 1.000 & 1.001   & 0.001 \\
	8&0.999 & 0.999 & 0.999 & 0.999   & 0.000 \\ \hline
\end{array}\]
It appears that the acceleration is rather high. To analyze the ultimate behavior of the error, we generalize the treatment in \cite[7.6]{foxparker}. Note that the function $\varphi_1(\theta_1,\theta_2)$ takes real values if $\theta_1=\theta_2$. Let $\alpha$ be a complex number such that $-2=1/\lambda_1=\varphi_1 (\alpha,\alpha)$. In our case, we can pick 
$\alpha \approx 1.925+\pi i$. Expressing (\ref{main}) in terms of the function $\varphi_1$, and focusing on the dominant terms, we see that this iterative process is ultimately equivalent to the following 
	\begin{align*}\label{main-v2}		
			\mb{y}^{(m)}=\ &(1+e^{-\alpha}+e^{-2\alpha})\left( M\mb{y}^{(m-1)} + \mb{g} \right) \\
			& -(e^{-\alpha}+e^{-2\alpha}+e^{-3\alpha})\left( \widetilde{M}\mb{y}^{(m-2)} + \widetilde{\mb{g}} \right)\\
			& +e^{-3\alpha}\mb{y}^{(m-3)}, \quad m\geq 3,
	\end{align*}
Let us set $a=1+e^{-\alpha}+e^{-2\alpha}$, $b=-(e^{-\alpha}+e^{-2\alpha}+e^{-3\alpha})$, and $c=e^{-3\alpha}$ for the further steps. Note that the error vector satisfies the iteration 
	\begin{equation*}\label{main-v4}
			\eta^{(m)} = aM\eta^{(m-1)} + b\widetilde{M}\eta^{(m-2)} +  c\eta^{(m-3)},
	\end{equation*}
which can expressed in the following matrix form
		\[\left[\begin{array}{c}
			\eta^{(m)}\\\eta^{(m-1)}\\\eta^{(m-2)}
		\end{array}\right] = \left[\begin{array}{ccc}
			aM& b\widetilde{M}&c\\\mb{I}&\mb{0}&\mb{0} \\ \mb{0}&\mb{I}&\mb{0}
		\end{array}\right] \left[\begin{array}{c}
			\eta^{(m-1)}\\ \eta^{(m-2)}\\ \eta^{(m-3)}
		\end{array}\right].\]
The ultimate convergence depends on the spectral radius of the matrix on the right-hand side. It is easy to see that the eigenvalues $\mu$ of this matrix satisfy the equation
		\[\mu^3- a\lambda\mu^2  -b\bar{\lambda}\mu - c=0,\]
where $\lambda$ is an eigenvalue of $M$. Recall that, in Example~\ref{ex1}, the largest $\mu$ occurred for the eigenvalue $\lambda$ whose absolute value is the spectral radius. We have a similar result in this example, as well. For $\lambda=\lambda_1$, we obtain $$|\mu|_\text{max} \approx 0.149.$$ 
This ratio exhibits a big improvement compared to the acceleration obtained in Example~\ref{ex1} with the classical semi-iterative Chebyshev method. 
		
We have verified experimentally that the ratio $\lVert\eta^{(m+1)}\rVert / \lVert\eta^{(m)}\rVert$ approaches to the number $0.149$ rather steadily. Even though the matrix $M$ is non-Hermitian, we may possibly have a strong relation between the spectral radius of $p_m(M)$ and $\lVert p_m(M) \rVert$, unlike the classical semi-iterative Chebyshev method. We believe that this is something interesting to study.
\end{example}

\section{Acknowledgement}
The author thanks Ö. Küçüksakallı for suggesting this problem. His support has accompanied the author throughout the preparation of this manuscript.

	{\small
		\def\refname{References}
		\newcommand{\etalchar}[1]{$^{#1}$}

	\end{document}